# Dual Darboux Frame of a Timelike Ruled Surface and Darboux Approach to Mannheim Offsets of Timelike Ruled Surfaces

**Mehmet Önder[1], H. Hüseyin Uğurlu[2]**

[1]*Celal Bayar University, Faculty of Arts and Sciences, Department of Mathematics, Muradiye Campus, Muradiye, Manisa, Turkey.* E-mail: mehmet.onder@cbu.edu.tr

[2]*Gazi University, Faculty of Education, Department of Secondary Education Science and Mathematics Teaching, Mathematics Teaching Program, Ankara, Turkey.* E-mail: hugurlu@gazi.edu.tr

**Abstract**
In this paper, we introduce the dual geodesic trihedron (dual Darboux frame) of a timelike ruled surface. By the aid of the E. Study Mapping, we consider timelike ruled surfaces as dual hyperbolic spherical curves and define the Mannheim offsets of timelike ruled surfaces by means of dual Darboux frame. We obtain the relationships between invariants of Mannheim timelike surface offsets. Furthermore, we give the conditions for these surface offsets to be developable.

## Introduction

In the surface theory, it is well-known that a surface is said to be "ruled" if it is generated by a continuously moving of a straight line in the space. Ruled surfaces are one of the simplest objects in geometric modeling that explains why these surfaces are one of the most fascinating topics of the surface theory and also used in many areas of science such as Computer Aided Geometric Design (CAGD), mathematical physics, moving geometry, kinematics for modeling the problems and model-based manufacturing of mechanical products. For example, the building materials such as wood are straight and they can be considered as straight lines. The result is that if engineers are planning to construct something with curvature, they can use a ruled surface since all the lines are straight [7].

An offset surface is offset a specified distance from the original along the parent surface's normal. Offsetting of curves and surfaces is one of the most important geometric operations in CAD/CAM due to its immediate applications in geometric modeling, numerical control (NC) machining, and robot navigation [4]. Especially, the offsets of the ruled surfaces have an important role in (CAGD) [12,13,14]. In [14], Ravani and Ku gave a generalization of the theory of Bertrand curve for Bertrand trajectory ruled surfaces on the line geometry. By considering the E. Study mapping, Küçük and Gürsoy have studied the integral invariants of closed Bertrand trajectory ruled surfaces in dual space [5]. They have given some characterizations of Bertrand offsets of trajectory ruled surfaces in terms of integral invariants (such as the angle of pitch and the pitch) of closed trajectory ruled surfaces.

Recently, a new offset of ruled surfaces has been defined by Orbay, Kasap and Aydemir [7]. They have called this new surface offset as Mannheim offset and obtained that every developable ruled surface has a developable Mannheim offset if and only if an equation should be satisfied between the geodesic curvature and the arc-length of spherical indicatrix of the reference surface. The corresponding characterizations of Mannheim offsets of ruled surfaces in Minkowski 3-space have been given in references [8,11]. Furthermore, in [9] Mannheim offsets of ruled surfaces have been studied in dual space with Blaschke approach and the characterizations of Mannheim offsets of ruled surfaces in terms of integral invariants of closed trajectory ruled surfaces have been introduced.

In this paper, we define the dual geodesic trihedron (dual Darboux frame) of a timelike ruled surface and give the real and dual curvatures of this surface. Then, we examine the Mannheim offsets of timelike ruled surfaces in view of their dual Darboux frame. Using dual representations of timelike ruled surfaces, we give some theorems and new results which characterize developable Mannheim timelike surface offsets.

**Preliminaries**

Let $\mathbb{R}_1^3$ be a 3-dimensional Minkowski space over the field of real numbers $\mathbb{R}$ with the Lorentzian inner product $\langle\,,\,\rangle$ given by

$$\langle \vec{a},\vec{a}\rangle = -a_1b_1 + a_2b_2 + a_3b_3,$$

where $\vec{a}=(a_1,a_2,a_3)$ and $\vec{b}=(b_1,b_2,b_3)\in \mathbb{R}^3$. A vector $\vec{a}=(a_1,a_2,a_3)$ of $\mathbb{R}_1^3$ is said to be timelike if $\langle \vec{a},\vec{a}\rangle < 0$, spacelike if $\langle \vec{a},\vec{a}\rangle > 0$ or $\vec{a}=0$, and lightlike (null) if $\langle \vec{a},\vec{a}\rangle = 0$ and $\vec{a}\neq 0$. Similarly, an arbitrary curve $\vec{\alpha}(s)$ in $\mathbb{R}_1^3$ is spacelike, timelike or lightlike (null), if all of its velocity vectors $\vec{\alpha}'(s)$ are spacelike, timelike or lightlike (null), respectively [6]. The norm of a vector $\vec{a}$ is defined by $\|\vec{a}\| = \sqrt{|\langle \vec{a},\vec{a}\rangle|}$. Now, let $\vec{a}=(a_1,a_2,a_3)$ and $\vec{b}=(b_1,b_2,b_3)$ be two vectors in $\mathbb{R}_1^3$. Then the Lorentzian cross product of $\vec{a}$ and $\vec{b}$ is given by

$$\vec{a}\times\vec{b}=(a_2b_3-a_3b_2,\ a_1b_3-a_3b_1,\ a_2b_1-a_1b_2).$$

By using this definition it can be easily shown that $\langle \vec{a}\times\vec{b},\vec{c}\rangle = -\det(\vec{a},\vec{b},\vec{c})$ [16,17].

The sets of the unit timelike and unit spacelike vectors are called hyperbolic unit sphere and Lorentzian unit sphere, respectively, and denoted by

$$H_0^2 = \{\vec{a}=(a_1,a_2,a_3)\in \mathbb{R}_1^3 : \langle \vec{a},\vec{a}\rangle = -1\},$$

and

$$S_1^2 = \{\vec{a}=(a_1,a_2,a_3)\in \mathbb{R}_1^3 : \langle \vec{a},\vec{a}\rangle = 1\},$$

respectively [16].

Let $D = \mathbb{R}\times\mathbb{R} = \{\bar{a}=(a,a^*): a,a^*\in \mathbb{R}\}$ be the set of pairs $(a,a^*)$. For $\bar{a}=(a,a^*)$, $\bar{b}=(b,b^*)\in D$ the following operations are defined on $D$:

Equality   : $\bar{a}=\bar{b} \Leftrightarrow a=b,\ a^*=b^*$
Addition   : $\bar{a}+\bar{b}=(a+b,\ a^*+b^*)$
Multiplication: $\bar{a}\bar{b}=(ab,\ ab^*+a^*b)$

Then the element $\varepsilon = (0,1)\in D$ satisfies following relationships

$$\varepsilon \neq 0,\ \varepsilon^2 = 0,\ \varepsilon 1 = 1\varepsilon = \varepsilon. \tag{1}$$

Let consider the element $\bar{a}\in D$ of the form $\bar{a}=(a,0)$. Then the mapping $f:D\to\mathbb{R},\ f(a,0)=a$ is an isomorphism. So, we can write $a=(a,0)$. By the multiplication rule we have that

$$\begin{aligned}
\overline{a} &= (a, a^*) \\
&= (a, 0) + (0, a^*) \\
&= (a, 0) + (0, 1)(a^*, 0) \\
&= a + \varepsilon a^*
\end{aligned}$$

Then $\overline{a} = a + \varepsilon a^*$ is called dual number and $\varepsilon$ is called dual unit. Thus the set of dual numbers is given by

$$D = \{\overline{a} = a + \varepsilon a^* : a, a^* \in \mathbb{R}, \varepsilon^2 = 0\}. \tag{2}$$

The set $D$ forms a commutative group under addition. The associative laws hold for multiplication. Dual numbers are distributive and form a ring over the real number field [3].

Dual function of dual number presents a mapping of a dual numbers space on itself. Properties of dual functions were thoroughly investigated by Dimentberg [2]. He derived the general expression for dual analytic (differentiable) function as follows

$$f(\overline{x}) = f(x + \varepsilon x^*) = f(x) + \varepsilon x^* f'(x), \tag{3}$$

where $f'(x)$ is derivative of $f(x)$ and $x, x^* \in \mathbb{R}$. This definition allows us to write the dual forms of some well-known functions as follows

$$\begin{cases} \cosh(\overline{x}) = \cosh(x + \varepsilon x^*) = \cosh(x) + \varepsilon x^* \sinh(x), \\ \sinh(\overline{x}) = \sinh(x + \varepsilon x^*) = \sinh(x) + \varepsilon x^* \cosh(x). \end{cases} \tag{4}$$

Let $D^3 = D \times D \times D$ be the set of all triples of dual numbers, i.e.,

$$D^3 = \{\tilde{a} = (\overline{a}_1, \overline{a}_2, \overline{a}_3) : \overline{a}_i \in D, i = 1, 2, 3\}. \tag{5}$$

Then the set $D^3$ is called dual space. The elements of $D^3$ are called dual vectors [1,3]. Similar to the dual numbers, a dual vector $\tilde{a}$ may be expressed in the form $\tilde{a} = \vec{a} + \varepsilon \vec{a}^* = (\vec{a}, \vec{a}^*)$, where $\vec{a}$ and $\vec{a}^*$ are the vectors of real space $\mathbb{R}^3$. Then for any vectors $\tilde{a} = \vec{a} + \varepsilon \vec{a}^*$ and $\tilde{b} = \vec{b} + \varepsilon \vec{b}^*$ of $D^3$, scalar product and cross product are defined by

$$\langle \tilde{a}, \tilde{b} \rangle = \langle \vec{a}, \vec{b} \rangle + \varepsilon \left( \langle \vec{a}, \vec{b}^* \rangle + \langle \vec{a}^*, \vec{b} \rangle \right), \tag{6}$$

and

$$\tilde{a} \times \tilde{b} = \vec{a} \times \vec{b} + \varepsilon \left( \vec{a} \times \vec{b}^* + \vec{a}^* \times \vec{b} \right), \tag{7}$$

respectively, where $\langle \vec{a}, \vec{b} \rangle$ and $\vec{a} \times \vec{b}$ are inner product and cross product of the vectors $\vec{a}$ and $\vec{b}$ in $\mathbb{R}^3$, respectively.

The norm of a dual vector $\tilde{a}$ is given by

$$\|\tilde{a}\| = \|\vec{a}\| + \varepsilon \frac{\langle \vec{a}, \vec{a}^* \rangle}{\|\vec{a}\|}, \quad (\vec{a} \neq 0). \tag{8}$$

A dual vector $\tilde{a}$ with the norm $1 + \varepsilon 0$ is called dual unit vector. The set of dual unit vectors is given by

$$\tilde{S}^2 = \{\tilde{a} = (\overline{a}_1, \overline{a}_2, \overline{a}_3) \in D^3 : \langle \tilde{a}, \tilde{a} \rangle = 1 + \varepsilon 0\}, \tag{9}$$

and called dual unit sphere (For details [1,3,18]).

The Lorentzian inner product of two dual vectors $\tilde{a} = \vec{a} + \varepsilon \vec{a}^*$, $\tilde{b} = \vec{b} + \varepsilon \vec{b}^* \in D^3$ is defined by

$$\langle \tilde{a}, \tilde{b} \rangle = \langle \vec{a}, \vec{b} \rangle + \varepsilon \left( \langle \vec{a}, \vec{b}^* \rangle + \langle \vec{a}^*, \vec{b} \rangle \right),$$

where $\langle \vec{a}, \vec{b} \rangle$ is the Lorentzian inner product of the vectors $\vec{a}$ and $\vec{b}$ in the Minkowski 3-space $\mathbb{R}_1^3$. Then a dual vector $\tilde{a} = \vec{a} + \varepsilon \vec{a}^*$ is said to be timelike if $\vec{a}$ is timelike, spacelike if $\vec{a}$ is spacelike or $\vec{a} = 0$ and lightlike (null) if $\vec{a}$ is lightlike (null) and $\vec{a} \neq 0$ [15].

The set of all dual Lorentzian vectors is called dual Lorentzian space and it is defined by
$$D_1^3 = \{ \tilde{a} = \vec{a} + \varepsilon \vec{a}^* : \vec{a}, \vec{a}^* \in \mathbb{R}_1^3 \}.$$

The Lorentzian cross product of dual vectors $\tilde{a}, \tilde{b} \in D_1^3$ is defined by
$$\tilde{a} \times \tilde{b} = \vec{a} \times \vec{b} + \varepsilon (\vec{a} \times \vec{b}^* + \vec{a}^* \times \vec{b}),$$
where $\vec{a} \times \vec{b}$ is the Lorentzian cross product in $\mathbb{R}_1^3$.

Let $\tilde{a} = \vec{a} + \varepsilon \vec{a}^* \in D_1^3$. Then $\tilde{a}$ is said to be dual timelike (resp. spacelike) unit vector if the vectors $\vec{a}$ and $\vec{a}^*$ satisfy the following equations:
$$<\vec{a}, \vec{a}> = -1 \ (resp. \ <\vec{a}, \vec{a}> = 1), \ <\vec{a}, \vec{a}^*> = 0. \tag{10}$$
The set of all dual timelike unit vectors is called the dual hyperbolic unit sphere, and is denoted by $\tilde{H}_0^2$,
$$\tilde{H}_0^2 = \{ \tilde{a} = (\bar{a}_1, \bar{a}_2, \bar{a}_3) \in D_1^3 : \langle \tilde{a}, \tilde{a} \rangle = -1 + \varepsilon 0 \}. \tag{11}$$
Similarly, the set of all dual spacelike unit vectors is called the dual Lorentzian unit sphere, and is denoted by $\tilde{S}_1^2$,
$$\tilde{S}_1^2 = \{ \tilde{a} = (\bar{a}_1, \bar{a}_2, \bar{a}_3) \in D_1^3 : \langle \tilde{a}, \tilde{a} \rangle = 1 + \varepsilon 0 \}. \tag{12}$$
(For details see [15]).

**Definition 1.** ([19]) *(i) Dual Hyperbolic angle:* Let $\tilde{x}$ and $\tilde{y}$ be dual timelike vectors in $D_1^3$. Then the dual angle between $\tilde{x}$ and $\tilde{y}$ is defined by $<\tilde{x}, \tilde{y}> = -\|\tilde{x}\|\|\tilde{y}\| \cosh \bar{\theta}$ and the dual number $\bar{\theta} = \theta + \varepsilon \theta^*$ is called the *dual hyperbolic angle*.

*(ii) Dual Lorentzian timelike angle:* Let $\tilde{x}$ be a dual spacelike vector and $\tilde{y}$ be a dual timelike vector in $D_1^3$. Then the dual angle between $\tilde{x}$ and $\tilde{y}$ is defined by $<\tilde{x}, \tilde{y}> = \|\tilde{x}\|\|\tilde{y}\| \sinh \bar{\theta}$ and the dual number $\bar{\theta} = \theta + \varepsilon \theta^*$ is called the *dual Lorentzian timelike angle*.

**Dual Representation and Dual Darboux Frame of a Timelike Ruled Surface**

In the Minkowski 3-space $\mathbb{R}_1^3$, an oriented timelike line $L$ is determined by a point $p \in L$ and a unit timelike vector $\vec{a}$. Then, one can define $\vec{a}^* = \vec{p} \times \vec{a}$ which is called moment vector. The value of $\vec{a}^*$ does not depend on the point $p$, because any other point $q$ in $L$ can be given by $\vec{q} = \vec{p} + \lambda \vec{a}$ and then $\vec{a}^* = \vec{p} \times \vec{a} = \vec{q} \times \vec{a}$. Reciprocally, when such a pair $(\vec{a}, \vec{a}^*)$ is given, one recovers the timelike line $L$ as $L = \{ (\vec{a} \times \vec{a}^*) + \lambda \vec{a} : \vec{a}, \vec{a}^* \in \mathbb{R}_1^3, \lambda \in \mathbb{R} \}$, written in parametric equations. The vectors $\vec{a}$ and $\vec{a}^*$ are not independent of one another and they satisfy the following relationships
$$\langle \vec{a}, \vec{a} \rangle = -1, \ \langle \vec{a}, \vec{a}^* \rangle = 0. \tag{13}$$

The components $a_i$, $a_i^*$ ($1 \leq i \leq 3$) of the vectors $\vec{a}$ and $\vec{a}^*$ are called the normalized Plucker coordinates of the timelike line $L$. From Eqs. (10), (11) and (13) we see that the dual timelike unit vector $\tilde{a} = \vec{a} + \varepsilon \vec{a}^*$ corresponds to timelike line $L$. This correspondence is known as E. Study Mapping: There exists a one-to-one correspondence between the timelike vectors of dual hyperbolic unit sphere $\tilde{H}_0^2$ and the directed timelike lines of the Minkowski space $\mathbb{R}_1^3$ [15]. By the aid of this correspondence, the properties of the spatial motion of a timelike line can be derived. Hence, the geometry of timelike ruled surfaces is represented by the geometry of dual hyperbolic curves lying on the dual hyperbolic unit sphere $\tilde{H}_0^2$.

Let now $(\tilde{e})$ be a dual hyperbolic curve represented by the dual timelike unit vector $\tilde{e}(u) = \vec{e}(u) + \varepsilon \vec{e}^*(u)$. The real unit vector $\vec{e}$ draws a curve on the real hyperbolic unit sphere $H_0^2$ and is called the (real) indicatrix of $(\tilde{e})$. We suppose throughout that it is not a single point. We take the parameter $u$ as the arc-length parameter $s$ of the real indicatrix and denote the differentiation with respect to $s$ by primes. Then we have $\langle \vec{e}', \vec{e}' \rangle = 1$. The vector $\vec{e}' = \vec{t}$ is the unit vector parallel to the tangent of the indicatrix. The equation $\vec{e}^*(s) = \vec{p}(s) \times \vec{e}(s)$ has infinity of solutions for the function $\vec{p}(s)$. If we take $\vec{p}_o(s)$ as a solution, the set of all solutions is given by $\vec{p}(s) = \vec{p}_o(s) + \lambda(s)\vec{e}(s)$, where $\lambda$ is a real scalar function of $s$. Therefore we have $\langle \vec{p}', \vec{e}' \rangle = \langle \vec{p}_o', \vec{e}' \rangle + \lambda$. By taking $\lambda = \lambda_o = -\langle \vec{p}_o', \vec{e}' \rangle$ we see that $\vec{p}_o(s) + \lambda_o(s)\vec{e}(s) = \vec{c}(s)$ is the unique solution for $\vec{p}(s)$ with $\langle \vec{c}', \vec{e}' \rangle = 0$. Then, the given dual curve $(\tilde{e})$ corresponding to the timelike ruled surface

$$\varphi_e = \vec{c}(s) + v\vec{e}(s), \tag{14}$$

may be represented by

$$\tilde{e}(s) = \vec{e} + \varepsilon \vec{c} \times \vec{e}, \tag{15}$$

where

$$\langle \vec{e}, \vec{e} \rangle = -1, \quad \langle \vec{e}', \vec{e}' \rangle = 1, \quad \langle \vec{c}', \vec{e}' \rangle = 0,$$

and $\vec{c}$ is the position vector of the striction curve. Then we have

$$\|\tilde{e}'\| = 1 - \varepsilon \det(\vec{c}', \vec{e}, \vec{t}) = 1 - \varepsilon \Delta, \tag{16}$$

where $\Delta = \det(\vec{c}', \vec{e}, \vec{t})$ which characterizes developable timelike surface, i.e, the timelike surface is developable if and only if $\Delta = 0$ [10]. Then, the dual arc-length $\bar{s}$ of the dual curve $(\tilde{e})$ is given by

$$\bar{s} = \int_0^s \|\tilde{e}'(u)\| du = \int_0^s (1 - \varepsilon \Delta) du = s - \varepsilon \int_0^s \Delta du. \tag{17}$$

From Eq. (17) we have $\bar{s}' = \dfrac{d\bar{s}}{ds} = 1 - \varepsilon \Delta$. Therefore, the dual unit tangent to the dual curve $(\tilde{e})$ is given by

$$\frac{d\tilde{e}}{d\bar{s}} = \frac{\tilde{e}'}{\bar{s}'} = \frac{\tilde{e}'}{1 - \varepsilon \Delta} = \tilde{t} = \vec{t} + \varepsilon(\vec{c} \times \vec{t}). \tag{18}$$

Introducing dual unit vector $\tilde{g} = -\tilde{e} \times \tilde{t} = \vec{g} + \varepsilon \vec{c} \times \vec{g}$, we have the dual frame $\{\tilde{e}, \tilde{t}, \tilde{g}\}$ which is known as dual geodesic trihedron or dual Darboux frame of $\varphi_e$ (or $(\tilde{e})$). Also, it is well known that the real orthonormal frame $\{\vec{e}, \vec{t}, \vec{g}\}$ is called the geodesic trihedron of the indicatrix $\vec{e}(s)$ with the derivations

$$\vec{e}' = \vec{t}, \quad \vec{t}' = \vec{e} + \gamma \vec{g}, \quad \vec{g}' = -\gamma \vec{t}, \tag{19}$$

where $\gamma$ is called the conical curvature [10].

Let now consider the derivatives of vectors of dual trihedron and find the dual Darboux formulae of a timelike ruled surface.

From Eq. (18) we have $\langle \tilde{t}, \tilde{t} \rangle = 1 + \varepsilon 0$. By using this equality and considering that $\tilde{g} = -\tilde{e} \times \tilde{t}$, we have

$$\left\langle \tilde{t}, \frac{d\tilde{t}}{ds} \right\rangle = 0, \quad \frac{d\tilde{g}}{ds} = -\tilde{e} \times \frac{d\tilde{t}}{ds}. \tag{20}$$

For the derivative of $\tilde{t}$ let write

$$\frac{d\tilde{t}}{ds} = \bar{\alpha}\tilde{e} + \bar{\beta}\tilde{t} + \bar{\gamma}\tilde{g}, \tag{21}$$

where $\bar{\alpha}, \bar{\beta}, \bar{\gamma}$ are dual functions of dual arc-length $\bar{s}$. The first equation of Eq. (20) gives that $\bar{\beta} = 0$. Thus from the second equation of Eq. (20) we have

$$\frac{d\tilde{g}}{ds} = -\bar{\gamma}\tilde{t}. \tag{22}$$

Finally, from Eq. (22) and the equality $\tilde{t} = \tilde{g} \times \tilde{e}$, we obtain

$$\frac{d\tilde{t}}{ds} = \tilde{e} + \bar{\gamma}\tilde{g}. \tag{23}$$

Then from Eqs. (18), (22) and (23) we have the following theorem:

**Theorem 1.** *The derivatives of the vectors of dual frame $\{\tilde{e}, \tilde{t}, \tilde{g}\}$ of a timelike ruled surface are given as follows,*

$$\frac{d\tilde{e}}{ds} = \tilde{t}, \quad \frac{d\tilde{t}}{ds} = \tilde{e} + \bar{\gamma}\tilde{g}, \quad \frac{d\tilde{g}}{ds} = -\bar{\gamma}\tilde{t}, \tag{24}$$

*and called dual Darboux formulae of timelike ruled surface $\tilde{e}$ (or $\varphi_e$). Then the dual Darboux vector of the trihedron is $\tilde{d} = -\bar{\gamma}\tilde{e} - \tilde{g}$.*

Let now give the invariants of the surface. Since $\bar{s}' = \frac{d\bar{s}}{ds} = 1 - \varepsilon\Delta$, from Eq. (22) it follows that

$$\tilde{g}' = -\bar{\gamma}(1 - \varepsilon\Delta)\tilde{t}. \tag{25}$$

On the other hand using the equality $\tilde{g} = \vec{g} + \varepsilon\vec{c} \times \vec{g}$, from Eq. (18) we have

$$\begin{aligned}\tilde{g}' &= -\gamma\vec{t} + \varepsilon(\vec{c}' \times \vec{g} - \gamma\vec{c} \times \vec{t}) \\ &= -\gamma\vec{t} - \varepsilon\gamma(\vec{c} \times \vec{t}) + \varepsilon(\vec{c}' \times \vec{g}) \\ &= -\gamma\tilde{t} + \varepsilon(\vec{c}' \times \vec{g})\end{aligned} \tag{26}$$

Therefore, from Eqs. (25) and (26) we obtain

$$-\bar{\gamma}(1 - \varepsilon\Delta)\tilde{t} = -\gamma\tilde{t} + \varepsilon(\vec{c}' \times \vec{g}), \tag{27}$$

which gives us

$$\bar{\gamma}(1 - \varepsilon\Delta) = \gamma + \varepsilon\delta, \tag{28}$$

where $\delta = \langle \vec{c}', \vec{e} \rangle$. Then from Eq. (28) we have

$$\bar{\gamma} = \gamma + \varepsilon(\delta + \gamma\Delta). \tag{29}$$

Moreover, since both $\vec{c}'$ and $\vec{e}$ are perpendicular to $\vec{t}$, for the real scalar $\mu$ we may write $\vec{c}' \times \vec{e} = \mu\vec{t}$. Then,

$$\Delta = \det(\vec{c}', \vec{e}, \vec{t}) = -\langle \vec{c}' \times \vec{e}, \vec{t} \rangle = -\mu \langle \vec{t}, \vec{t} \rangle = -\mu.$$

Hence we have $\vec{e} \times (\vec{c}' \times \vec{e}) = -\Delta \vec{e} \times \vec{t} = \Delta \vec{g}$ and $\vec{c}' = -\delta \vec{e} + \Delta \vec{g}$.

The functions $\gamma(s)$, $\delta(s)$ and $\Delta(s)$ are the invariants of the timelike ruled surface $\varphi_e$. They determine the timelike ruled surface uniquely up to its position in the space. For example, if $\delta = \Delta = 0$ we have $\vec{c}$ is constant. It means that the timelike ruled surface $\varphi_e$ is a timelike cone.

**Elements of Curvature of a Dual Hyperbolic Curve**

Similar to the Euclidean geometry, in the Lorentzian geometry, the radius of curvature, $R$, of a curve at a point is a measure of the radius of the Lorentzian circular arc which best approximates the curve at that point. It is the inverse of the curvature. Then if we consider the dual Lorentzian curve, the dual radius of curvature can be obtained by considering dual versions of well-known formulae which gives radius of curvature of a timelike or a spacelike curve. Then the dual radius of curvature of dual hyperbolic curve (timelike ruled surface) $\tilde{e}(s)$ is can be calculated analogous to common Lorentzian differential geometry of curves as follows

$$\bar{R} = \frac{\left\| \dfrac{d\tilde{e}}{ds} \right\|^3}{\left\| \dfrac{d\tilde{e}}{ds} \times \dfrac{d^2\tilde{e}}{ds^2} \right\|} = \frac{1}{\sqrt{|1 - \bar{\gamma}^2|}}. \tag{30}$$

The unit vector $\tilde{d}_o$ with the same sense as the Darboux vector $\tilde{d} = -\bar{\gamma}\tilde{e} - \tilde{g}$ is given by

$$\tilde{d}_o = -\frac{\bar{\gamma}}{\sqrt{|1 - \bar{\gamma}^2|}} \tilde{e} - \frac{1}{\sqrt{|1 - \bar{\gamma}^2|}} \tilde{g}. \tag{31}$$

It is clear that $\tilde{d}_o$ is timelike (or spacelike) if $|\bar{\gamma}| > 1$ (or $|\bar{\gamma}| < 1$). Then, the dual angle between $\tilde{d}_o$ and $\tilde{e}$ satisfies the followings,

$$\begin{cases} \cosh \bar{\rho} = -\dfrac{\bar{\gamma}}{\sqrt{|1 - \bar{\gamma}^2|}}, \quad \sinh \bar{\rho} = -\dfrac{1}{\sqrt{|1 - \bar{\gamma}^2|}}, \text{ if } |\bar{\gamma}| > 1, \\ \sinh \bar{\rho} = -\dfrac{\bar{\gamma}}{\sqrt{|1 - \bar{\gamma}^2|}}, \quad \cosh \bar{\rho} = -\dfrac{1}{\sqrt{|1 - \bar{\gamma}^2|}}, \text{ if } |\bar{\gamma}| < 1. \end{cases} \tag{32}$$

where $\bar{\rho}$ is the dual spherical radius of curvature. Hence,

$$\bar{R} = \begin{cases} -\sinh \bar{\rho}, \text{ if } |\bar{\gamma}| > 1, \\ -\cosh \bar{\rho}, \text{ if } |\bar{\gamma}| < 1, \end{cases} \text{ and } \bar{\gamma} = \begin{cases} \coth \bar{\rho}, \text{ if } |\bar{\gamma}| > 1, \\ \tanh \bar{\rho}, \text{ if } |\bar{\gamma}| < 1. \end{cases}$$

**Darboux Approach to Mannheim Offsets of Timelike Ruled Surfaces**

Let $\varphi_e$ be a timelike ruled surface generated by dual timelike unit vector $\tilde{e}$ and $\varphi_{e_1}$ be a ruled surface generated by dual unit vector $\tilde{e}_1$ and let $\{\tilde{e}(\bar{s}), \tilde{t}(\bar{s}), \tilde{g}(\bar{s})\}$ and $\{\tilde{e}_1(\bar{s}_1), \tilde{t}_1(\bar{s}_1), \tilde{g}_1(\bar{s}_1)\}$ denote the dual Darboux frames of $\varphi_e$ and $\varphi_{e_1}$, respectively. Then, $\varphi_e$ and $\varphi_{e_1}$ are called Mannheim surface offsets, if

$$\tilde{g}(\bar{s}) = \tilde{t}_1(\bar{s}_1) \tag{33}$$

holds, where $\bar{s}$ and $\bar{s}_1$ are the dual arc-lengths of $\varphi_e$ and $\varphi_{e_1}$, respectively. By this definition, it is seen that the Mannheim offset of $\varphi_e$ is also a timelike ruled surface, but the generator of this surface can be timelike or spacelike. In this study, we consider the Mannheim offset surface with timelike ruling. If one considers the second situation, ruling is spacelike, and uses the Definition 1 (ii) similar results can be found.

Let now $\varphi_e$ and $\varphi_{e_1}$ be timelike ruled surfaces with timelike rulings and let these surfaces form a Mannheim surface offset. Then the relation between the trihedrons of the ruled surfaces $\varphi_e$ and $\varphi_{e_1}$ can be given as follows

$$\begin{pmatrix} \tilde{e}_1 \\ \tilde{t}_1 \\ \tilde{g}_1 \end{pmatrix} = \begin{pmatrix} \cosh\bar{\theta} & \sinh\bar{\theta} & 0 \\ 0 & 0 & 1 \\ \sinh\bar{\theta} & \cosh\bar{\theta} & 0 \end{pmatrix} \begin{pmatrix} \tilde{e} \\ \tilde{t} \\ \tilde{g} \end{pmatrix}, \tag{34}$$

where $\bar{\theta} = \theta + \varepsilon\theta^*$, $(\theta, \theta^* \in \mathbb{R})$ is the dual hyperbolic angle between dual timelike generators $\tilde{e}$ and $\tilde{e}_1$ of $\varphi_e$ and $\varphi_{e_1}$. The angle $\theta$ is called the offset angle which is the angle between the rulings $\vec{e}$ and $\vec{e}_1$, and $\theta^*$ is called the offset distance which is measured from the striction point $\vec{c}$ of $\varphi_e$ to the striction point $\vec{c}_1$ of $\varphi_{e_1}$. From Eq. (33) we write $\vec{c}_1 = \vec{c} + \theta^*\vec{g}$. Then, $\bar{\theta} = \theta + \varepsilon\theta^*$ is called dual hyperbolic offset angle of the Mannheim timelike ruled surfaces $\varphi_e$ and $\varphi_{e_1}$. If $\theta = 0$, then the Mannheim timelike surface offsets are said to be oriented offsets.

Now, we give some theorems and results characterizing Mannheim timelike surface offsets. When we write $\varphi_e$ and $\varphi_{e_1}$, we mean timelike ruled surfaces with timelike ruling and for short we don't write the Lorentzian characters of the surfaces hereinafter.

**Theorem 2.** *Let $\varphi_e$ and $\varphi_{e_1}$ form a Mannheim surface offset. The offset angle $\theta$ and the offset distance $\theta^*$ are given by*

$$\theta = -s + c, \quad \theta^* = \int_0^s \Delta du + c^*, \tag{35}$$

*respectively, where $c$ and $c^*$ are real constants.*

**Proof.** Suppose that $\varphi_e$ and $\varphi_{e_1}$ form a Mannheim offset. Then from Eq. (34) we have

$$\tilde{e}_1 = \cosh\bar{\theta}\tilde{e} + \sinh\bar{\theta}\tilde{t}. \tag{36}$$

By differentiating Eq. (36) with respect to $\bar{s}$ we have

$$\frac{d\tilde{e}_1}{d\bar{s}} = \sinh\bar{\theta}\left(1 + \frac{d\bar{\theta}}{d\bar{s}}\right)\tilde{e} + \cosh\bar{\theta}\left(1 + \frac{d\bar{\theta}}{d\bar{s}}\right)\tilde{t} + \bar{\gamma}\sinh\bar{\theta}\tilde{g}. \tag{37}$$

Since $\frac{d\tilde{e}_1}{d\bar{s}}$ and $\tilde{g}$ are linearly dependent, by using Theorem 1 and Eq (33), from Eq. (37) we get $\frac{d\bar{\theta}}{d\bar{s}} = -1$. Then for the dual constant $\bar{c} = c + \varepsilon c^*$ we write $\bar{\theta} = -\bar{s} + \bar{c}$ and from Eq. (17) we have

$$\theta = -s + c, \quad \theta^* = \int_0^s \Delta du + c^*,$$

where $c$ and $c^*$ are real constants.

From Eq. (35), the following corollary can be given.

**Corollary 1.** *Let $\varphi_e$ and $\varphi_{e_1}$ form a Mannheim surface offset. Then $\varphi_e$ is developable if and only if $\theta^* = c^* = \text{constant}$.*

**Theorem 3.** *Let $\varphi_e$ and $\varphi_{e_1}$ form a Mannheim surface offset. Then there is the following differential relationship between the dual arc-length parameters of $\varphi_e$ and $\varphi_{e_1}$*

$$\frac{d\overline{s}_1}{d\overline{s}} = \overline{\gamma} \sinh \overline{\theta}. \tag{38}$$

**Proof.** Suppose that $\varphi_e$ and $\varphi_{e_1}$ form a Mannheim offset. Then, from Theorem 1 and Theorem 2 we have

$$\frac{d\tilde{e}_1}{d\overline{s}_1} = \tilde{t}_1 = \overline{\gamma} \sinh \overline{\theta} \frac{d\overline{s}}{d\overline{s}_1} \tilde{g}. \tag{39}$$

From Eq. (33) we have $\tilde{t}_1 = \tilde{g}$. Then Eq. (39) gives us

$$\overline{\gamma} \sinh \overline{\theta} \frac{d\overline{s}}{d\overline{s}_1} = 1, \tag{40}$$

and from Eq. (40) we get Eq. (38).

**Corollary 2.** *Let $\varphi_e$ and $\varphi_{e_1}$ form a Mannheim surface offset. Then there are the following relationships between the real arc-length parameters of $\varphi_e$ and $\varphi_{e_1}$*

$$\frac{ds_1}{ds} = \gamma \sinh \theta, \quad \frac{ds ds_1^* - ds^* ds_1}{ds^2} = \theta^* \gamma \cosh \theta + (\delta + \gamma \Delta) \sinh \theta. \tag{41}$$

**Proof.** Let $\varphi_e$ and $\varphi_{e_1}$ form a Mannheim surface offset. Then by Theorem 3, Eq. (38) holds. By considering Eqs. (4) and (29), the real and dual parts of Eq. (38) are

$$\frac{ds_1}{ds} = \gamma \sinh \theta, \quad \frac{ds ds_1^* - ds^* ds_1}{ds^2} = \theta^* \gamma \cosh \theta + (\delta + \gamma \Delta) \sinh \theta. \tag{42}$$

In Corollary 1, we give the relationship between the offset distance $\theta^*$ and developable timelike ruled surface $\varphi_e$. Now we give the condition for $\varphi_{e_1}$ to be developable according to $\theta^*$. From Eq. (17) we have

$$ds^* = -\Delta ds, \quad ds_1^* = -\Delta_1 ds_1. \tag{43}$$

Then writing Eq. (43) in Eq. (42) and using Eq. (41) we get

$$\Delta_1 = -\left(\theta^* \coth \theta + \frac{\delta}{\gamma}\right), \tag{44}$$

which gives the following corollaries:

**Corollary 3.** *Let $\varphi_e$ and $\varphi_{e_1}$ form a Mannheim surface offset. Then Eq. (44) holds.*

**Corollary 4.** *Let $\varphi_e$ and $\varphi_{e_1}$ form a Mannheim surface offset. Then $\varphi_{e_1}$ is developable if and only if $\theta^* = -\dfrac{\delta}{\gamma}\tanh\theta$ holds.*

**Theorem 4.** *Let $\varphi_e$ and $\varphi_{e_1}$ form a Mannheim surface offset. There exists the following relationship between the invariants of the surfaces and offset angle $\theta$, offset distance $\theta^*$,*

$$\delta_1 = \frac{\delta}{\gamma}\coth\theta - \theta^*. \tag{45}$$

**Proof.** Let the striction lines of $\varphi_e$ and $\varphi_{e_1}$ be $c(s)$ and $c_1(s_1)$, respectively, and let $\varphi_e$ and $\varphi_{e_1}$ form a Mannheim surface offset. Then, from the Mannheim condition we can write

$$\vec{c}_1 = \vec{c} + \theta^* \vec{g}. \tag{46}$$

Differentiating Eq. (46) with respect to $s_1$ we have

$$\frac{d\vec{c}_1}{ds_1} = \left(\frac{d\vec{c}}{ds} - \theta^*\gamma\vec{t} + \frac{d\theta^*}{ds}\vec{g}\right)\frac{ds}{ds_1}. \tag{47}$$

We know that $\delta_1 = \langle d\vec{c}_1/ds_1, e_1\rangle$. Then from Eqs. (34) and (47) we obtain

$$\delta_1 = \left(\cosh\theta\langle d\vec{c}/ds, \vec{e}\rangle + \sinh\theta\langle d\vec{c}/ds, \vec{t}\rangle - \theta^*\gamma\sinh\langle\vec{t},\vec{t}\rangle\right)\frac{ds}{ds_1}. \tag{48}$$

Since $\delta = \langle d\vec{c}/ds, e\rangle$ and $\langle d\vec{c}/ds, \vec{t}\rangle = 0$, from Eq. (48) we write

$$\delta_1 = (\delta\cosh\theta - \theta^*\gamma\sinh\theta)\frac{ds}{ds_1}. \tag{49}$$

Furthermore, from Eq. (41) we have

$$\frac{ds}{ds_1} = \frac{1}{\gamma\sinh\theta}, \tag{50}$$

Substituting Eq. (50) in Eq. (49) we obtain

$$\delta_1 = \frac{\delta}{\gamma}\coth\theta - \theta^*.$$

**Theorem 5.** *If $\varphi_e$ and $\varphi_{e_1}$ form a Mannheim surface offset, then the relationship between the conical curvature $\gamma_1$ of $\varphi_{e_1}$ and offset angle $\theta$ is given by*

$$\gamma_1 = -\coth\theta. \tag{54}$$

**Proof.** From Eqs. (19) and (34) we have

$$\begin{aligned}\gamma_1 &= -\langle\vec{g}_1', \vec{t}_1\rangle \\ &= -\left\langle\frac{d}{ds_1}(\sinh\theta\vec{e} + \cosh\theta\vec{t}), \vec{g}\right\rangle \\ &= -\gamma\cosh\theta\frac{ds}{ds_1}\end{aligned} \tag{55}$$

From the first equality of Eqs. (41) and (55) we have $\gamma_1 = -\coth\theta$.

**Theorem 6.** *If the surfaces $\varphi_e$ and $\varphi_{e_1}$ form a Mannheim surface offset, then the dual conical curvature $\bar{\gamma}_1$ of $\varphi_{e_1}$ is obtained as follows*

$$\bar{\gamma}_1 = -\coth\theta + \varepsilon\left(\frac{2\delta}{\gamma}\coth\theta + \theta^*\operatorname{cosech}^2\theta\right). \tag{56}$$

**Proof.** From Eqs. (29), (44), (45) and (55) by direct calculation we have Eq. (56).

**Theorem 7.** *If the surfaces $\varphi_e$ and $\varphi_{e_1}$ form a Mannheim surface offset, then the dual curvature of $\varphi_{e_1}$ is given by*

$$\bar{R}_1 = \sinh\theta - \varepsilon\sinh\theta\cosh\theta\left(\frac{2\delta}{\gamma}\cosh\theta + \theta^*\operatorname{cosech}\theta\right). \tag{57}$$

**Proof.** From Eq. (56) we have

$$\sqrt{\left|1 - \bar{\gamma}_1^2\right|} = \operatorname{cosech}\theta + \varepsilon\coth\theta\left(\frac{2\delta}{\gamma}\cosh\theta + \theta^*\operatorname{cosech}\theta\right).$$

Then from Eq. (30) we have Eq. (57).

**Theorem 8.** *Let $\varphi_e$ and $\varphi_{e_1}$ form a Mannheim surface offset and let $|\bar{\gamma}_1| > 1$. Then, the dual spherical radius of curvature of $\varphi_{e_1}$ is given by*

$$\cos\bar{\rho}_1 = \cosh\theta - \varepsilon\sinh^2\theta\left(\frac{2\delta}{\gamma}\cosh\theta + \theta^*\operatorname{cosech}\theta\right). \tag{58}$$

**Proof.** From (32) we have $\sinh\bar{\rho}_1 = -\bar{R}_1$. Then considering Eqs. (4) and (57) we get

$$\sinh\rho_1 = -\sinh\theta, \quad \rho_1^*\cosh\rho_1 = \cosh\theta\sinh\theta\left(\frac{2\delta}{\gamma}\cosh\theta + \theta^*\operatorname{cosech}\theta\right). \tag{59}$$

By using the trigonometric relationship $\cosh^2\rho_1 - \sinh^2\rho_1 = 1$, from the first equality of Eq. (59) we have

$$\cosh\rho_1 = \cosh\theta. \tag{60}$$

Then by writing Eq. (60) in the second equality of Eq. (59) we obtain

$$\rho_1^* = \sinh\theta\left(\frac{2\delta}{\gamma}\cosh\theta + \theta^*\operatorname{cosech}\theta\right). \tag{61}$$

By the aid of the extension $\cosh\bar{\rho}_1 = \cosh\rho_1 + \varepsilon\rho_1^*\sinh\rho_1$, from Eqs. (59)-(61) we have Eq. (58).

From Eq. (60) we have the following corollary:

**Corollary 5.** *Let $\varphi_e$ and $\varphi_{e_1}$ form a Mannheim surface offset. Then the offset angle $\theta$ is equal to the real spherical radius of curvature $\rho_1$, i.e., $\rho_1 = \theta$.*

If we assume that $|\bar{\gamma}_1| < 1$, then we have equalities for a timelike ruled surface whose Darboux vector is spacelike and the obtained equalities will be analogue to given ones.

**Conclusions**

The dual geodesic trihedron (dual Darboux frame) of a timelike ruled surface is introduced. Then the characterizations of Mannheim offsets of timelike ruled surfaces are given in view of dual Darboux frame and new relations between the invariants of Mannheim timelike surface offsets are obtained. The relationships for Mannheim timelike surface offsets to be developable according to offset angle and offset distance are given. The results of the paper are new characterizations of Mannheim timelike surface offsets.

**References**

[1] Blaschke W (1945) Differential Geometrie and Geometrischke Grundlagen ven Einsteins Relativitasttheorie Dover. New York
[2] Dimentberg F.M (1965) The Screw Calculus and its Applications in Mechanics. English translation: AD680993, Clearinghouse for Federal and Scientific Technical Information, (Izdat. Nauka, Moscow, USSR)
[3] Hacısalihoğlu H.H (1983) Hareket Geometrisi ve Kuaterniyonlar Teorisi. Gazi Üniversitesi Fen-Edb Fakültesi
[4] Hoschek J, Lasser D (1993) Fundamentals of computer aided geometric design. Wellesley MA:AK Peters
[5] Küçük A, Gürsoy O (2004) On the invariants of Bertrand trajectory surface offsets. App Math and Comp 151:763-773
[6] O'Neill B (1983) Semi-Riemannian Geometry with Applications to Relativity. Academic Press London
[7] Orbay K, Kasap E, Aydemir I (2009) Mannheim Offsets of Ruled Surfaces. Mathematical Problems in Engineering Article ID 160917
[8] Önder M, Uğurlu H.H, On the Developable of Mannheim offsets of timelike ruled surfaces in Minkowski 3-space. arXiv:0906.2077v5 [math.DG]
[9] Önder M, Uğurlu H.H, Some Results and Characterizations for Mannheim Offsets of the Ruled Surfaces. arXiv:1005.2570v3 [math.DG]
[10] Önder M, Uğurlu H.H (2012) Frenet Frames and Invariants of Timelike Ruled Surfaces. Ain Shams Eng J http://dx.doi.org/10.1016/j.asej.2012.10.003
[11] Önder M, Uğurlu H.H, Kazaz M, Mannheim offsets of spacelike ruled surfaces in Minkowski 3-space. arXiv:0906.4660v3 [math.DG]
[12] Papaioannou S.G, Kiritsis D, (1985) An application of Bertrand curves and surfaces to CAD/CAM. Computer Aided Design 17(8):348-352
[13] Pottmann H, Lü W, Ravani B (1996) Rational ruled surfaces and their offsets. Graphical Models and Image Processing 58(6):544-552
[14] Ravani B, Ku T.S (1991) Bertrand Offsets of ruled and developable surfaces. Comp Aided Geom Design 23(2):145-152
[15] Uğurlu H.H, Çalışkan A (1996) The Study Mapping for Directed Spacelike and Timelike Lines in Minkowski 3-Space $\mathbb{R}_1^3$. Mathematical and Computational Applications 1(2):142-148
[16] Uğurlu H.H, Çalışkan A (2012) Darboux Ani Dönme Vektörleri ile Spacelike ve Timelike Yüzeyler Geometrisi. Celal Bayar Üniversitesi Yayınları Yayın No: 0006
[17] Uğurlu H.H (1999) The relations among instantaneous velocities of trihedrons depending on a spacelike ruled surface. Hadronic Journal 22:145-155
[18] Veldkamp G.R (1976) On the use of dual numbers, vectors and matrices in instantaneous spatial kinematics. Mech and Mach Theory 11:141-156
[19] Yakut N.N (2012 Reel ve Dual Uzaylarda Açı Kavramı. CBÜ Fen Bilimleri Enstitüsü, Yüksek Lisans Tezi